\documentclass[12pt]{amsart}
 
\usepackage{amsfonts, amsbsy, amsmath, amsthm, amssymb, latexsym, enumerate,verbatim}

\newtheorem{propo}{Proposition}
\newtheorem{defi}[propo]{Definition}

\newtheorem{lemma}[propo]{Lemma}
\newtheorem{corol}[propo]{Corollary}

\newtheorem{theo}[propo]{Theorem}

\newtheorem{prob}{Problem}

\newcommand{\bpr}{\begin{prob}\label}
\newcommand{\epr}{\end{prob}}

\newcommand{\bl}{\begin{lemma}\label}
\newcommand{\el}{\end{lemma}}

\newcommand{\ra}{ \rightarrow }
\newcommand{\lan}{ \langle }
\newcommand{\ran}{ \rangle }

\newcommand{\ld}{, \ldots ,}

\newcommand{\CC}{\mathbb C} 
\newcommand{\FF}{\mathbb F}
\newcommand{\ZZ}{\mathbb Z} 

\newcommand{\End}{\mathop{\rm End}\nolimits}
\newcommand{\Hom}{\mathop{\rm Hom}\nolimits}

\newcommand{\Id}{\mathop{\rm Id}\nolimits}

\newcommand{\Irr}{\mathop{\rm Irr}\nolimits}

\newcommand{\al}{\alpha}

\newcommand{\ep}{\varepsilon}

\newcommand{\lam}{\lambda }
\newcommand{\si}{\sigma }

\newcommand{\om}{\omega }
\newcommand{\Om}{\Omega }

\newcommand{\GG}{\mathbf{G}}

\newcommand{\up}{^{-1}}
\def\d12{{_{12}}}

\def\ac{{algebraically closed }}
\def\acf{{algebraically closed field }}

\def\au{{automorphism }}

\def\cc{{characteristic }}

\def\ei{{eigenvalue }}
\def\eis{{eigenvalues }}

\def\f{{following }}

\def\ii{{if and only if~\,}}

\def\ir{{irreducible }}
\def\irt{{irreducible. }}
\def\irr{{irreducible representation }}

\def\itf{{It follows that }}

\def\mult{{multiplicity }}

\def\po{{polynomial }}
\def\pos{{polynomials }}
\def\pot{{polynomial. }}

\def\rep{{representation }}
\def\reps{{representations }}
\def\rept{{representation. }}
\def\repst{{representations. }}

\def\ag{algebraic group }

\newcommand{\med}{\medskip}

\newcommand{\bp}{\begin{proof} }
\newcommand{\enp}{\end{proof}}

\date{}

\begin{document}

 \title[Fixed point free actions of group elements]
{Fixed point free actions of group elements
   in linear groups and finite group representations  }

\bigskip

\author{A.E. Zalesski}

\address{Department of Mathematics,  University of Brasilia, Brasilia-DF, Brazil}
\email{alexandre.zalesski@gmail.com}




\begin{abstract} This is an expository paper aimed to outline the current situation with problems related with the occurrence of \ei $1$ of elements in linear groups and group representations. 
\end{abstract}

\maketitle



\bigskip
\centerline{Dedicated to the memory of Martin Isaacs}


\bigskip

\tableofcontents

\section{Introduction}

In the theory of transformation groups a significant role belongs to the study of invariants, that is, the objects remaining unchanged under the group action. In group \reps the invariants are the non-zero underlying space vectors fixed by every group element. A more specific problem is the study of invariants for individual group elements in group representations, that is, invariants in this case are the non-zero elements of the 1-eigenspace of an element in question. Many concrete questions came from applications,  whereas a systematic study is still at
 the infant level I think.

We restrict our discussion   to the elements $g\in GL_n(F)$, $F$ is a field, that are diagonalizable under a field extension. They are called semisimple and, if $g$ is of finite order $|g|$,  are exactly those with $|g|$ coprime to the characteristic of $F$.  Properties of non-semisimple linear group elements is another important area of intensive study, which we do not discuss in this survey.  Note that this article can be viewed as a continuation of author's expository paper \cite{z09}. In this survey we outline recent achievements and new areas of study.

The area of studing \eis is very wide and I start from general overview. Currently 
there are several projects which aim on understanding regularities related with 
eigenvalues  of elements of finite linear groups and 
group \repst Ideally, given a finite linear group $G\subset GL_n(F)$, 
one would wish to obtain full information on all eigenvalues of all elements $g\in G$.
This is possible for many groups of bounded order, for instance, those listed in the
 Atlas of Brauer characters \cite{bat} of finite groups. Some more groups can be dealt with
using computer programs available in famous systems GAP, MAGMA and CHEVIE.  
This approach is limited by computer powers. Obtaining regularities at the first stage
requires specifying situations where strong regularities exist, and stating certain conjectures.

In this survey I try to outline areas of intensive research and give short comments
on main achievements. I do not definitely pretend to be complete.

I usually assume that $F$ is algebraically closed; this is not a significant assumption as the general case to large extent reduces to this one. The case $F$ of characteristic 0 is much easier and regularities are much stronger. However, I mostly assume $F$ to be of arbitrary characteristic $\ell$.

In the most general context the problem of our discussion can be stated as follows:

 \begin{prob}\label{gg99}  Given an algebraically closed field $F$, a linear group $G\subset GL_n(F)$ and an element $g\in G$, determine  the Jordan normal form of $g$, in particular, 
the  \eis of $g$ and their multiplicities.
 \end{prob}


As mentioned above, there may be many ways to deal with Problem \ref{gg99} if the triple $(F,n,G,g)$ is explicitly given. However, we are interested with regularities which hold for many groups
and large classes of their elements. To reveil regularities, even at the level of a conjecture, one needs to specify   the entries $(F,n,G,g)$ or/and limit oneself with a less precise answer. 
For instance, one can ask for determining the degree of the minimal polynomial of $g$ or occurrence of \ei 1  for every  element of $G$ or for specific class of them. 

The most important restriction of general nature is that of irreducibility of $G$. This will be held in almost all cases considered below. Other natural restrictions considered in the literature are:

\med
$(i)$ $n$ is small or of specific nature (for instant a prime);

$(ii)$ $g$ does not lie in a proper normal subgroup of $G$, equivalently, $G$ is generated by the elements conjugate to $g;$

$(iii)$  $g$ is a $p$-element for a prime $p$, or even $|g|=p;$

$(iv)$ $G$ or $G'$ is quasisimple, sometimes specified to be an alternating or classical group or a group of Lie type;

$(v)$ the characteristic $\ell$ of $F$ does not divide $|G|;$  if $G$ is of Lie type one differs the cases where $\ell$ equals or distinct from the defining characteristic of $G;$

$(vi)$ $g$ has a simple spectrum, that is, all \eis of $g$ are of \mult 1; more generally, at most one \ei  of $g$ is of \mult greater than $1;$  

$(vii)$ the degree $\deg g$ of the minimal \po of $g$ is less than $|g|$ (or  $o(g)$, the degree of $g$ in $G/Z(G));$ 

$(viii)$ $g$ does not have \ei 1. 

\med
There are ¨uniform¨ (in a sense) versions of the above restrictions, for instance, one asks when  all elements of $G$ have \ei 1.  

One is also interested with lower and upper bounds for the \ei multiplicities of certain elements. Similar questions are discussed for group representations.

In this paper I focus on problems and results that concern with \ei 1 properties of $g\in G\subset GL_n(F)$ in notation of Problem \ref{gg99}, especially with the existence of \ei 1.

\med
{\it  
Notation} $F$ is a field of \cc $\ell$, $F^\times$ is the multiplicative group of $F$; this is assumed to be algebraically closed unless otherwise stated. $GL_n(F)$ is the group of $(n\times n)$-matrices over $F$, and we often write $GL_n(F)=GL(V)$ in order to say that $V$ is the underlying vector space for $GL_n(F)$.
Notation for other classical groups are standard, we follow those in \cite{KL}, as well as 
for simple groups. $\FF_q$ is a finite field of $q$ elements. 

If $G$ is a group then $|G|$ is the order of $G$, $Z(G)$ is the center and $G'$ the derived subgroup of $G$. The identity element of $G$ is usually denoted by $1$. If $g,h\in G$ then $[g,h]:=ghg\up h\up$. By $o(g)$ we denote the order of $g\in G$ in $G/Z(G)$. If $p$ is a prime then $p'$-elements of $G$ are those of order coprime to $p$. If $S\subset G$ is a subset then $C_G(S)=\{g\in G: [x,s]=1\}$ and $N_G(S)=\{g\in G: gS g\up=S\}$. If $p$ is a prime   then $O_p(G)$ is the maximal normal $p$-subgroup of $G$. Notation for groups of Lie type are standard, we mainly follow those of \cite{atl}. The symmetric and alternating groups on $n$ points are denoted by $S_n$ and $A_n$, respectively. 

If $g\in GL_n(F)=GL(V)$ then   $\deg g$ denotes the degree of the mi-nimal \po of $g$.
We say that $g$ is fixed point free if 1 is not an \ei of $g$. The  1-eigenspace of $g$
is often introduced as $C_V(g)$.   An element $g$ with \ei 1 is sometimes called unisingular.

Let $\rho:G\ra GL_n(F)$ be a \rept We often write $\rho\in \Irr_F(G)$ to say that $\phi$ is \irt The trivial \rep of $G$ is denoted by $1_G$, and we write $\rho_G^{reg}$ for the regular \rep of $G$. If $H$ is a subgroup of $G$, we write $\phi|_H$ for the restriction of $\phi$ to $H$. If $\lam $ is a \rep of $H$, we denote by $\lam^H$ the induced \rept

\section{Motivation and some applications}

In this section we briefly outline some areas of potential applications of results concerning properties of \ei 1 of linear groups or group \repst 

\subsection{Linear algebra and geometry}

Formally, the \ei study is a part of linear algebra. However, regularities of the \ei behavior in 
linear groups are mostly of very different nature. The facts that a linear group element $g$
does not have \ei 0, and the \eis of $g$ are $|g|$-roots of unity, are practically insignificant. Note that the study of \eis of an element $g\in GL_n(F)=GL(V)$ is equivalent to the study of cyclic group $\lan g \ran$, more precisely, the study of \ir constituents of $\lan g \ran$, and their multiplicities.

This hints that one cannot expect strong regularities without assuming that    
a group $G \subset GL_n(F)$ in question is \irt So we assume below that elements 
of interest belong to an \ir subgroup $G\subset GL_n(F)$. Moreover, we assume that
$G$ is a proper subgroup, and  does not contain $SL_n(F)$. In addition, 
the assumption that $G$ is a classical group is a very modest restriction on elements $g\in G$. Therefore, one  usually excludes the cases where  $G$   
contains the derived subgroup of a classical group on $V$.  

The linear algebra view point suggests some natural questions on elements $g\in G\subset GL_n(F)$ to be considered:

\med
Question 1. What can be said on the minimal \pos of elements  $g\in G,$ in particular, on $\deg g$? 

\med
Question 2. What can be said on the minimal and maximal \mult of the \eis of $g$?

\med
Question 3. When is 1  an \ei of $g?$ When does it occur with \mult 1?

\med
Question 4. When are all \ei multiplicities of $g$ equal to 1?  
 
\med
Question 5. When is a certain  $|g|$-root of unity an \ei of $g$? 

\med
These questions specify some situations where well described regularities can be expected. There are variations of these questions that specify the classes of groups $G$ (simple or solvable, say) or   particular types of elements (of prime order, say). In the opposite direction
one can turn to a more uniform version asking when some \ei property holds for all group elements or for a certain class of them. The following two problems illustrate this approach.
 
\begin{prob}\label{11}  For an \ir finite simple subgroup $G\subset GL_n(F)$ and a prime $p$ dividing $|G|$, when is the minimal \po degree $\deg g$ of an element $g\in G$ of order $p$ less than $p$?\end{prob}

Note that if 1 is not an \ei of $g$ then $\deg g<p$. Therefore,  one can single out
the \f special case:

\begin{prob}\label{12}  For an \ir finite simple subgroup $G\subset GL_n(F)$ and a prime $p$ dividing $|G|$, when is $1$ not an \ei of an element $g\in G$ of order $p$?\end{prob}

Let $G\subset GL_n(F)$ be a finite \ir group. Obviously,  $\deg g\leq o(g)$. In general,   the
situation where   $\deg g= o(g)$ is common, whereas  $\deg g<o(g)$ is exceptional.
The  problem of describing all exceptions is hardly treatable. Problem \ref{11}
is still open, many special cases are resolved in \cite{TZ08} and \cite{TZ22}. Probably, 
the \f special case is one of most difficult:

\bpr{13} Let $p$ be a prime, $G=A_p$, the alternating group, and let $g\in G$ be of order p. Let $F$ be an \acf of characteristic $\ell\neq 0,p$. Determine the \ir $F$-\reps $\phi$ of $G$ such that   $\deg\phi(p)<p$. 
\epr

We expect that $\dim\phi<p$. Note that the case where $\ell>n/2$ follows from a result by Thompson \cite[Lemma 3.2]{Th1}.   In this case Sylow $\ell$-subgroups of $G$ are cyclic. 
Note that the general theory of \reps for groups with cyclic $\ell$-subgroup has  powerful tools to study the \ei problems.


As we show below,  the following vector space covering problem (stated explicitly in \cite{cer})  is partially connected with the \ei 1 problem: 

\bpr{is19} 
Let $H\subset GL(V)$ be a subgroup and $W$ a subspace of $V$ over a finite field $\FF_q$.  When is $V=\cup_{h\in H} h(W)$?\epr

More generally, given a normal subgroup $N$ of a group $G$, one  asks 
when $N=\cup_{g\in G}\, gKg\up$ for a proper subgroup $K$ of $N$.  If $N=G$ then $N$  never coincides with the union of the conjugates of $K$. This fact goes back to C. Jordan (1872).

\subsection{Arithmetic algebraic geometry}

Motivated by applications to a certain aspect of Galois theory and abelian varieties Katz \cite[Problem II, p. 483]{k81} stated a linear group problem which we reformulate as follows:

\bpr{p4}   Let $p$ be a prime and  $H=Sp_{2n}(p)$  the symplectic group.
 Describe \ir subgroups G of H whose every element has \ei $1$.  \epr

\med
Specifying $H$ to be a symplectic group over a prime field is not natural from group theory point of view. So we state a general problem which can highlight other special cases:

\bpr{4pm}  Describe  \ir  finite subgroups $G\subset GL_n(F)$ whose all elements have \ei $1$.\epr

We call such $G$ unisingular.  For small $n$ Problems \ref{p4} and \ref{4pm} can be solved by inspection of finite \ir subgroups of $GL_n(F)$, see \cite{cu10,cu19,cz} and \cite{z24} and the bibliography there. For more comments on  applications to algebraic geometry see \cite{cu09,cu12,cu24}.
In general, Problem \ref{4pm} is not treatable as there are too many groups in question. One can try to reduce the search by considering maximal subgroups:  

\bpr{p5}  Describe maximal 
 \ir locally  finite subgroups $G\subset GL_n(F)$ whose all elements have \ei $1$.
\epr

Note that  in Problem \ref{p5} $G$ may be the union of an infinite series $G_1\subset G_2\subset \cdots  $ of \ir finite subgroups $G_i\subset GL_n(F),$ $i=1,2,...$. Groups $SO_3(\ell^{n_i})$, $\ell$ odd, 
and $n_i|n_{i+1}$ form in a sense one of the simplest examples. In any case $G$ is contained in a group isomorphic to $GL_n(\overline{F}_\ell)$, where $\overline{F}_\ell$ is the algebraic closure of $\FF_\ell$. 
If $G$ is infinite then one can consider the connected component $G^\circ$ of $G$ with respect of Zarisky topology (see \cite[\S 8]{Di});  
 this is a normal subgroup of finite index in $G$. Furthermore, one can consider  the  Zarisky topology closure $\overline{G}$ of $G$ in  $GL_n(\overline{F}_\ell)$. 
Note that $\overline{G}$ is unisingular if so is $G$ (as $g\in GL_n(F)$ has \ei 1 \ii
$g$ satisfies the polynomial equation $\det (x\cdot \Id-g)=0$). In addition, $\overline{G}$
is an \ag over  $\overline{F}_\ell$   and $\overline{G}^\circ$ is a connected algebraic group over 
$\overline{F}_\ell$ (which can be reducible). Nonetheless, Problem \ref{p5}
is closely connected with the following problem for algebraic group representations:

\bpr{p5a}  Let $\mathbf{G}$ be a simple algebraic group that is not finite. Determine 
the   rational \ir \reps of $\mathbf{G}$ whose all elements have \ei $1$.\epr

The above example is related to the fact that the simple algebraic group $SL_2(\overline{F}_\ell)$, $\ell>2$,  has a unisingular  \irr of degree 3.  Problem \ref{p5a} can be reduced
to some problem on the weight 0 existence in a given \irr of a simple algebraic group.    
More details are given in Section 5.1. 

Another version of Problem \ref{4pm} is the following:

\bpr{p5b}  Describe minimal  \ir   finite subgroups $G\subset GL_n(F)$ that contain a fixed point free element.
\epr
 
We pay readers' attention to the recent book \cite{KT} as a good source of linear group problems related to algebraic geometry.

\subsection{Game theory}

\begin{theo}\label{isb} {\rm (Isbell \cite{Is}, see also \cite[p. 10 ]{bg16})} Let $n>0$ be even. Then there exists an $n$-player homogeneous  game \ii there exists a transitive permutation group on $n$ points that  contains no fixed point free $2$-element (in other words, every $2$-element fixes a point).
\end{theo}

In fact, there is a bijective correspondence between $n$-player homogeneous  games
and transitive permutation groups on $n$ points that  contain no fixed point free $2$-element.  This naturally leads to a problem of determining the transitive permutation subgroups of $S_n$, $n$ even,  that  contain no fixed point free $2$-element.  Specifying $p=2$ is not natural from permutation group theory view point. So one can ask
(see \cite[Problem 5]{z24}):


\bpr{is1} Given  a natural number $n>1$ and a prime $p|n$ determine   transitive subgroups $G$ of $S_n$ whose every $p$-element fixes a point. Equivalently,  
determine   transitive subgroups $G$ of $S_n$ containing a fixed point free $p$-element.\epr

This problem is hardly treatable in full generality. Current experience is not sufficient in order to state any conjecture. Note that, for a given $n$ and a transitive subgroup $G\subset S_n$
there exists a fixed point free $p$-element  for {\it } some prime $p|n$, see  \cite{fks} or \cite[Theorem 1.3.1]{bg16}. So the feature of Problem \ref{is1} is that $p$ is fixed in advance. 
In \cite{Be} the author shortly mentions a related problem: given a transitive permutation group $G$ of $ S_n$ determine primes $p|n$ such that $g$ has a fixed point free $p$-element. Note that fixed point free elements in $S_n$ are often referred as derangements in literature. 

The existence of a derangement of $p$-power order does not always imply that of order $p$.  However, there are very few  primitive  transitive permutation group that have no derangement of prime order, see 
\cite{bow,bg16,gu3}.  

Linear groups over  finite fields containing a fixed point free $p$-element are  good sources of examples of groups $G$ for Problem \ref{is1}. Indeed, let $G\subset GL_n(q)=GL(V)$
and $g\in G$ a $p$-element. Suppose that $g$ does not have \ei 1 on  $V$. Then $g$
fixes no element of $\Om=V\setminus \{0\}$. Therefore, $g$ acts fixed point freely on every non-trivial $G$-orbit on $\Om$. This fact is exploited in \cite{np8}.  

Problem \ref{is1} looks  analogous to Problem \ref{12}, however, the nature of these problems are   somehow different. For instance, it is not true that a finite \ir group 
contains a fixed point free element. However, in some special cases these two problems are equivalent. For instance let $G$ be a minimal non-abelian solvable group that is not of prime power order. Then $G$ has an elementary abelian normal $r$-subgroup $A$ for some prime $r$ and $G/A$ is of prime order $p$, say. In addition, $A$ is the only non-trivial proper normal subgroup of $G$.

\begin{propo}\label{es3} {\rm See \cite[Lemma 3.9]{z24}} Let $G=AM$ be a semidirect product of an elementary abelian r-group $A$ and a  group M of prime order $ p\neq r$.  Suppose that $A$ has no non-trivial proper normal subgroup of $G$. 
Then the following are equivalent: 

$(1)$ $G$ is isomorphic to an \ir subgroup of $X\subset GL_d(r)=GL(V)$ for some d and every element of $X$ has \ei $1$; 

$(2 )$ $G$ is isomorphic to a transitive subgroup of $Y\subset S_{rp}=Sym(\Om)$ and every r-element of $Y$ fixes a point on $\Om;$

$(3)$ Let $K$ be a subgroup of index $r$ in $A$. Then $A=\cup_{g\in G}\, gKg\up$.  
\end{propo}

The equivalence of $(1)$, (2) and (3) explains     that Problems \ref{12}, \ref{is1} and \ref{is19} have some common roots.
Note that the group $Y$ in  item (2) of Proposition \ref{es3} is imprimitive, so methods developed in \cite{bow,bg16,gu3} do not work to decide whether $Y$ has a derangement of order $p$.
 Observe that the group $G$ depends on two parameters $r,p$.  See Section 7 for details.

\subsection{Group theory}

There are many aspects of studying \eis but those I shall discuss go back to Hall-Higman's theorem on the minimal polynomial degrees of elements in \reps of $p$-solvable groups. 

\begin{theo}\label{hh2} {\rm \cite[Theorem B]{hh}} Let $\ell=p>2$, let $G\subset GL_n(F)$ be a p-solvable group with no non-trivial normal p-subgroup and $g\in G$ a p-elements of order $p^a$. Then the minimal \po degree of $g$ is at least $(p-1)p^{a-1}$. \end{theo}

In this theorem $F$ is an arbitrary field of characteristic $p$, in particular, it is true for $|F|=p$. This result  plays a substantial role in finite group theory, and used in the proof of the classification theory of finite simple groups. So this result attracted a lot of attention, see \cite{Br,Cf,Fw}, \cite[Ch. IX, \S 2]{HB2} etc.

 If   $G\subset GL_n(F)$ is an \ir $p$-solvable group over an \acf $F$ of characteristic $\ell>0$ then $G$  lifts to characteristic 0, see \cite[Ch. X, Theorem 2.1]{Fe}.

The notion of lifting plays an essential role in the study of \ei problems of finite linear groups over fields of prime characteristic. One says that a finite subgroup $G\subset GL_n(F)$
is {\it  liftable} (or that $G$ lifts to characteristic 0) if 
there is an isomorphism $\eta:G\ra  GL_n(\CC)$ such the $\ell$-Brauer character of $G$ coincides on the $\ell'$-elements with the character of $\eta(G)$.  (Note that every finite  $\ell'$-subgroup of  $ GL_n(F)$ is liftable.)
So if $g\in G$ is   an $\ell'$-element   then the \ei multiplicities of   $g$ are the same as those of $\eta(g)$. 

In turn, if $G\subset GL_n(\CC) $ is a finite group then for every $\ell>0$ there exists a homomorphism $ \si:G\ra  GL_n(F)$ such that $\deg g\geq \deg\si(g)$ (a Brauer reduction of $G$ modulo $\ell$). Therefore, 
Theorem \ref{hh2} remains valid  for $G\subset  GL_n(\CC)$ with $O_p(G)=1$. 

 Therefore, for $p$-solvable groups Theorem \ref{hh2}
extends to \ir  $p$-solvable   group $G\subset GL_n(F)$ with $  O_p(G)=1$ over a field $F$ of arbitrary characteristic $\ell$. Indeed, $\deg g=\deg\eta(g)$ if $\ell\neq 0,p$ and 
$\deg\eta(g)\geq \si(\eta(g))$. 

Is the assumption of $G$ $p$-solvable substantial? At least one can state a general form problem:

\bpr{is15} {\rm (General Hall-Higman type problem)}  
 Given a finite \ir linear group   $G\subset GL_n(F)$,  determine   $p$-elements  $g\in G$ such that $\deg g<o(g)$, where $o(g)$ is the order of $g$ in $G/Z(G)$ and $\deg g$ is the degree of minimal \po  of $g$. \epr 

 For quasisimple groups $G$ the main source of results on Problem \ref{is15} is \cite{TZ08} and  \cite{TZ22}  together with  earlier works \cite{DZ1} and \cite{DM08}.

 Similarly, we can state:

\bpr{is11} {\rm  
(General minimal \po problem)} Given an \ir finite subgroup $G\subset GL_n(F)$, determine the minimal \po degrees of its elements. \epr
 
One can interpret Problem \ref{is11} in a stronger sense as determining the pairs $(g,\deg g)$, $g\in G$. This version is much harder.  

One assumes that the ground field $F$ and the degree $n$ is given.   It is clear that a full solution to
Problem \ref{is11} can only be obtained for some small values of $n$. So one has to select some special cases where there can exist a describable  solution. One can fix an upper bound for $n$, or single out a series of $n$ with
a strong restriction on $n$, for instance, for $n$ prime, or impose a restriction to the groups in question, say, $p$-groups or simple groups, or narrow attention to  elements 
  of particular nature such as reflections, of prime order etc. Much depends on the characteristic of $F$.   

Instead of climbing up on $n$, one can consider Problem \ref{is11} from the representation theory point of view:
 
\bpr{is12}
 Given a group $G$ and an \irr $\phi: G\ra GL_n(F)$ determine the minimal \po degrees of $\phi(g)$ for $g\in G$.\epr
 
Formally, Problems \ref{is11} and \ref{is12} are equivalent. However, one can  single out series of groups $G$  which can be treated uniformly. In addition, there is a lot of information on \ir \reps of individual groups, especially for those that are regarded as important. 

If $G\subset GL_n(F) $ and an element $g\in G$ does not have \ei 1 then $\deg g\leq |g|-1$.  As we wish to focus to  \ei 1 problems, we single out this special case:

  \begin{prob}\label{t91}   Given an \ir finite subgroup $G\subset SL_n(F)$,  determine the elements $g\in G$ that have \ei $1$. \epr

There is an evidence based on certain experience 
 that vast majority of elements of $G$ have \ei value 1, and those without \ei 1 are exceptional and can   in principle be determined. However, 
one has to exclude the cases where $G$ is a quasisimple group of Lie type in characteristic $\ell$ and the symmetric and alternating groups $S_m$ and $A_m$ with $m<n+3$. In this cases the situation is different, and one can turn to studying the proportion of the number of elements having \ei 1. 

The absence of \ei 1 means that the element in question acts fixed point freely on the underlying vector space (formally, the zero vector is the only one fixed by $g$). 
This property makes sense when one considers an action of one group $G$, say, on another group $H$  by automorphisms. If $g\in G$ and and $C_H(g)=1$, one says that $g$ acts fixed point freely. There are a significant literature devoted to fixed  point free  automorphisms  of solvable groups.

Problem \ref{t91} relates with a probability problem discussed in \cite{gt03}:

  \begin{prob}\label{t92}   Given an \ir finite subgroup $G\subset GL_n(F)$, compute the probability $P(G)$ for an element   $g\in G$ to have  \ei $1$. Obtain a lower bound for $P(G)$. \epr

Note that $P(G)$ is the ratio $e/|G|$, where $e$ is the number of elements $g\in G$ that have \ei 1. In fact, Problem  \ref{t91} is stated and discussed in \cite{gt03} for the case where $G$ is a      simple group of Lie type in defining characteristic $\ell$. Obviously, $P(G)\geq 1/|G|$; the groups $G$ with  $P(G)=1/|G|$ are determined in \cite{Zss}. 
If $\ell>0$ then one can be interested in determining subgroups  $G\subset GL_n(F)$ whose all $\ell'$-element do not have \ei 1. This problem has been solved in \cite{FLT}.

Another problem in general group theory which requires to study the existence of 
\ei 1 of certain group elements in  group representations arises in the study of 
the   Grunberg-Kegel graph $\Gamma(G)$ of a finite group $G$ (known also as the prime   graph of $G$). The vertices of  $\Gamma(G)$ are the primes dividing $|G|$, and  primes $p,r$ are connected by an edge \ii $G$ has an element of order $pr$.  Clearly,    $\Gamma(G)$ accumulates very small information of the structure of $G$. Surprizingly, in some non-trivial cases the graph determines the group in full. 

\smallskip
The following problem is discussed in numerios papers (see \cite{K,Ma} and references there):  

\bpr{40} Describe  groups $G$ that are determined by $\Gamma(G)$.
\epr  

Let $\om(G)$ denote the set of element orders of a finite group $G$. Then $\om(G)$ determines $\Gamma(G)$. There are a significant number of publications concerning with the problem of determining a group $G$ by  $\om(G)$   \cite{M,D} etc.  For an illustration of the connection of Problem \ref{40} with Problem \ref{is14}, consider the following situation.  
Let $H=G/A$, where $A$ is an elementary abelian $p$-group and $C_H(A)=A$. Suppose that $H$ has 
no element of order $pr$. Then the vertices $p,r$ are not connected on $\Gamma(H)$, and these are connected on $\Gamma(G)$ \ii some $p$-element $x\in G$ centralizes an non-identity  element of $A$. 
Viewing $A$ as a vector space over $\FF_p$ and $H$ as a subgroup of $GL(A)$, one interprets this by saying that $|C_A(x)|>1$ \ii $x$ has \ei 1 as an element of $GL(A)$.

\section{Fixed point free elements of linear groups and the occurrence of eigenvalues 1 of elements in group representations: Overview of research area}

For reader's convenience we recall the following definition.

\begin{defi}\label{du1} An element $g\in GL_n(F)$ is called unisingular if it has \ei $1$. A linear group  $G$ is called unisingular if all elements of $G$ have \ei $1$. 
A \rep $\phi$ of a group $G$ is called unisingular if $\phi(g)$ has \ei $1$ for every $g\in G$.\end{defi}

Note that Problem \ref{4pm} is in a sense  a uniform  version of Problem \ref{t91}, which is much easier and better treatable. Moreover, we single out the case of $G$ simple: 

\bpr{u5i} Determine unisingular  \ir finite simple subgroups  of $GL_n(F)$.    \epr

Guralnick and Tiep \cite{gt03} address to a  problem of even  more uniform nature:

\bpr{is25} 
Given a  field $F$ determine all finite simple  groups $G$ whose {\it all}
  $F$-\reps are unisingular.\epr

They obtained in \cite{gt03} a solution for $G$ of Lie type in defining characteristic $p>0$  and $F$ of characteristic $p$. In general, Problem \ref{is25} remains open. 

As explained in the previous section, the case of $G$ simple is basic; alternatively, one can single out concrete groups and classes of elements for analyzing. For instance one can restrict  the analysis to $p$-elements, see Problem \ref{is15}.  

To approach Problem \ref{u5i} one has to use the classification of finite simple groups and the \rep theory of every finite simple group. So we rewrite this problem as follows:  

  \begin{prob}\label{t41}  Let G be a simple group 
and F an algebraically closed field.  
Determine unisingular \ir \reps $\phi$ of G over F.
\end{prob} 

Not much  is  known on  Problem \ref{t41} if $F$ is the field of complex numbers. In \cite{ap2} and \cite{st1} this has been solved for alternating and symmetric groups.  
Some sufficient conditions of unisingularity for ordinary representations of finite groups of Lie type were given in \cite{Z90}. In particular, every ordinary \rep of simple groups $E_8(q),F_4(q),{}^2F_4(q)$ and $G_2(q)$ with $q$ odd  is unisingular. More such groups are given in \cite{gt03}, for instance $PSp_{2n}(q)$ for $q$ an odd prime.

We mention one more problem in spirit of Problem \ref{t41}. An element $g\in G$ is called {\it rational} if $g$ conjugate to $g^i$ whenever $(i,|g|)=1$, equivalently,  the cyclic groups $\lan g \ran$
and $\lan g^ i\ran$ coincide. If $g$ is  rational and $\lam$ is an \ei of $g$ then so is $\lam^i$ with  $(i,|g|)=1$. This hints that for rational elements the \ei 1 problem is much better treatable
than in general. A special case considered in  \cite{zv23} supports this view.

\bpr{is16} 
  Given a finite group  $G$ and a field $F$, determine the  \ir $F$-representations of   $G $  in which every rational element of $G$ has \ei $1$.    \epr

This is equivalent to determining the \ir \reps in which some rational element acts fixed point freely. Recall that all elements of the symmetric groups $S_n$ are rational, so for this group Problem \ref{is16} is equivalent to Problem \ref{t41}.

In most cases Problem \ref{t41} is expected to have a describable solution. 
It is unlikely however that a similar problem is treatable for arbitrary finite groups. Nonetheless we state it below to outline the framework which various special cases belong to.  

\bpr{is13}  Given a   finite group $G$ and a field $F$, determine the   \ir $F$-representations $\phi$ of  $G$  such that $1$ is  an \ei of $\phi(g)$ for {\it every} $g\in G$.  \epr

Much more is known on the occurrence of \ei 1  of $p$-elements $g\in G$ due to results on the Hall-Higman type problems, see \cite{TZ08}. 
We find it convenient to single out this case  as follows:

\bpr{is14}
Given a finite group   $ G$ and a field $F$, determine the   \ir $F$-representations of   $G$ in  which some $p$-element of $G$ does not  have \ei $1$.    \epr


  Problem \ref{is14} and  Problem \ref{is15} must have similar solutions unless $G$ has a normal subgroup isomorphic to a quasisimple group of Lie type in characteristic $\ell$.
In this case  Problem \ref{is15} looks untreatable, whereas  Problem \ref{is14} can have a describable solution in a few rather general cases (see \cite{z18},  \cite{z22} and \cite{TZ25}). 
 However, in any case one can expect rather  strong regularities under relatively light restrictions on $G$.   At this stage it looks reasonable to  turn to various special cases in order to increase the understanding of the general problem and build some machinery and tools to deal with more general situations.

\med
In order to understand the role of $n$ in Problem \ref{4pm},  one may consider the following problem:

\bpr{is21} Determine integers $n>1$ such that $GL_n(F)$
has no unisingular \ir  subgroup.  \epr

For $\ell\neq 2$ we have the following result.

\begin{theo}\label{th5}{\rm \cite[Theorem 1.2]{z24}} Suppose that   $\ell\neq 2$. Then $GL_n(F)$
has an absolutely \ir   unisingular subgroup unless $n=2,4$ or $\ell=3$ and $n=8$. \end{theo} 

So Problem \ref{is21} reduces to $p=2$, and this case looks rather complex.
 The case of  $p=2$ is discussed in several papers by J. Cullinan together with 
    a more  special problem going back to N. Katz \cite{k81}:

\bpr{is22} Determine unisingular   subgroups of $Sp_{2n}(2)$, the symplectic group over a field of two elements. \epr 

This problem is motivated by some applications to an arithmetic aspect of algebraic geometry, see Section 2.2.  
Cullinan himself solved Problem \ref{is22} for $n=2,3$; the case of $n=4$ has been completed in \cite{cz}.

 If $n$ is large, the number of groups giving a solution to  Problem \ref{is22}  dramatically increases  together with $n$, especially for $n$ a 2-power.   
 Problem \ref{is22} remains of great interest for particular values of $n$. 
One can singled out some asymptotic aspect of Problem \ref{is22}:

\bpr{is23}  
$ (a)$  Determine integers $n>1$ such that $Sp_{2n}(2)$ has no \ir  unisingular  subgroup. $(b)$ Are there infinitely many $n$ with this property?   \epr

One easily observes that if $G\subset Sp_{2n}(2)$ is a unisingular  subgroup then 
the wreath product $G \wr S_m$ is   a unisingular  subgroup of $Sp_{2mn}(2)$
for every integer $m>0$.
This observation leads to a more specific  

\bpr{is24}    Are there infinitely many {\it primes} $p$ such that $Sp_{2p}(2)$ has no unisingular \ir subgroup?  \epr

\med
 
The following result is obtained in \cite[Theorem 1.3]{z24}:

\begin{theo}\label{th6} Let $1<n<125$.

If $n\in\{1,2,3,5,9,27,29,43,53,81,89,113\}$ then $Sp_{2n}(2)$ has no  unisingular subgroup. 

If $n\notin\{1,2,3,5,9,27,29,43,53,81,113, 47,58,  67,83,86,103,107\}$ then 
$Sp_{2n}(2)$ contains a unisingular   subgroup.\end{theo}

\section{Some  observations on minimal \pos and \ei 1 of elements of linear groups}

If $g\in G\subset GL_n(F)$ then $\deg g$   cannot exceed the order $|g|$ of $g$. 
If $g$ is diagonalizable then $\deg g$ is the number of distinct \eis of $g$, and if 1
is not an \ei of $g$ then $\deg g<|g|$. One can easily verify the \f elementary observation. 

\bl{ssu} Let $g\in GL_n(F)$ be an element of finite order, and $g=su$, where $s,u\in \lan g \ran$, $u$ is an $\ell$- and $s$ is an $\ell'$-element. Then $g$ has \ei  $1$   \ii s has.  \el

 If $G$ is finite then the knowledge of its character table for $\ell=0$ allows to obtain 
full information on the \eis of $\phi(g)$ for $g\in G$. If $\ell>0$ and  the Brauer character table of $G$ is known then one can similarly obtain full information on the \eis of $\phi(g)$ for $\ell'$-elements of $g\in G$.
Irreducible finite quasisimple  subgroups of $GL_n(F)$ for $n\leq 250$ are listed  in \cite{HMa}.  

One can be interested in a more ambitious problem of determining full list of the \eis of elements of a group in question. However, since a scalar multiple of   a matrix has the same
 minimal polynomial degree but not the same list of eigenvalues,  Problems \ref{is11} and \ref{is12} have a more uniform solution, and hence better treatable than the problem of computing all eigenvalues.   According with our experience, knowledge of 
$\deg g$ allows to easily determine the \eis themselves by a more or less routine computations. 

 Clearly, any result depends on the conjugacy class of an element, and by induction reason it is natural to assume that this does not lie in a proper normal subgroup of the group in question. There are some configurations where no reasonable regularities can be expected. As mentioned in section 2.1,  this is the case for the general linear group $G=GL_n(F)$, where one can exhibit an element with almost arbitrary a priory possible minimal \pot (For instance, the minimal \po degree cannot exceed $n$. In addition, if $g$ is conjugate in $G$ to $g^m$ for some $m>1$ then 
$\lam^m$ is an \ei of $g$ whenever  $\lam$ is.) This concerns also with $SL_n(F_1)$ and, in a practical sense, the classical groups such as
$SO_n(F_1)$, $Sp_n(F_1)$ and $SU_n(F_1)$ for some subfield $F_1$ of $F$.
To some extent this also concerns with tensor products of these groups such that 
$GL_k(F)\otimes GL_l(F)$. Because of this, and also for the  induction purpose one often assumes
$G$ to be tensor-indecomposable, that is, $G$ does not lies in $GL_k(F)\otimes GL_l(F)$
for $k,l>1, kl=n$. There is a reason to exclude imprimitive groups. 

The  following lemma  reduces eigenvalue  problems to much more special situations:

\bl{hh1}  Let $V=V_1 \oplus \cdots \oplus V_k$ be a direct sum of vector spaces
over a field $F$, and let $g\in  GL(V)$ be an element of prime power order permuting
transitively the subspaces $V_1\ld  V_k$.
 Set $G =\lan g\ran$ and $G^k = \lan g^k\ran$. Then the following hold:

$(i)$ Assume $(g, \ell) = 1$. If $\mu$ is an eigenvalue of $g^k$ of multiplicity $s$, then every $k$-th-root of $\mu$ is an eigenvalue of g of multiplicity $s/k$. (Observe
that $s/k$ is also the multiplicity of $\mu$ on every $V_i$.)

$(ii)$ $\deg g = k \cdot  \deg g^k$. \el

This fact leads to an observation often referred as Higman's lemma (see \cite[Theorem 1.10]{HB2}:

\bl{h2h} Let $G\subset GL(V)$ be a subgroup and $g\in G$. Suppose that $|g|$ is a prime power, g normalizes an abelian subgroup $A$ of G, $(|A|,\ell)=1$ and  $C_{\lan g \ran}(A)=1$. 
Then $\deg g=|g|$. 
\el

In fact, one can drop  the assumption $C_{\lan g \ran}(A)=1$ and then replace the conclusion by $\deg g\geq k$, where $k$ is the order of $g$ modulo $C_G(A)$.   In order to deduce Lemma \ref{h2h}
from Lemma \ref{hh1} decompose $V$ as a direct sum of homogeneous $FA$-modules;  these   are permuted by $g$.  As $|g|$ is a prime power, there is an orbit of size $k$, and we can apply Lemma \ref{hh1} to the sum of the spaces of this orbit.  This reasoning explains in part why dealing with prime power order elements $g\in G$ is easier than in general; if $g\in S_n$ is not a prime power then $g$ may not have an orbit of size $|g |$. 
 
Lemmas \ref{hh1} and \ref{h2h}  are  useful on an initial stage of  analysis of Problem \ref{11} (and also Problem \ref{is12}). Say, if $|g|=p$ is a prime and $g\in G$ normalizes an abelian subgroup without centralizing it, we have $\deg g=p$, in particular, $g$ has \ei 1.  As an illustration, every $p$-element of $SL_n(p)$ for $n>2$ has this property. In a more general context, this leads  to 
the \f result. To state it,  for $p$ odd, we denote by $\Delta_ 1(p)$ (respectively: $\Delta_2(p)$) the set $ 1 \cup  \{\ep ^j \}$, where $1 \neq \ep\in F$,
$\ep^p = 1$ and $j$ runs over the non-squares (respectively: the squares) of $\ZZ$  modulo $p$.

Denote by ${\rm Spec}\, g$ the set of distinct \eis of an element $g\in GL_n(F)$.
Note that if $G$ contains a scalar matrix $z$ then $\deg g=\deg zg$.

\begin{propo}\label{z90}  {\rm  (Cf. \cite{Z86} and \cite{Z88}.)} Let H be a quasi-simple group of Lie type in defining characteristic p, such that $(p,|Z(H)|) = 1$, and let $g\in H$ be an element of order p. Let $\theta$ be a faithful
irreducible $\ell$-representation of H for $\ell\neq p$, and suppose that $\deg \theta (g) < p$. Then p is odd and one of the following holds:

$(1)$ $ H = PSU(3,p)$, $\dim\theta = p(p -1)$ and g is a transvection;

$(2)$ $H = SL(2,p^2)$, $\dim\theta = (p^2 -1)/2;$

$(3)$ $H = Sp(4,p),$ $\dim\theta = (p^2 - 1)/2$  and g is not a transvection;

$(4)$ $ H = PSp(4,p)$, $\dim\theta = p(p - 1)2/2$ and g is a transvection;

$(5)$ $ H = Sp(2n,p)$ or $PSp(2n,p)$, $n > 1$, $\dim\theta = (p^n \pm  1)/2$, g is a transvection and
${\rm Spec} \,\theta(g) = \Delta_1(p)$                         or $\Delta_2(p)$. If  $\ell = 2$, then $H = PSp(2n,p)$ and only the minus sign
has to be taken in the expression for $\dim\theta$. If  $\ell \neq  2$, then $H = Sp(2n,p)$ if $\dim\theta$ is
even, while $H = PSp(2n,p)$ if $\dim\theta$ is odd;

$(6)$ $H = SL(2,p)$ or $PSL(2,p)$, and either $ \dim\theta = (p + 1)/2 $ with ${\rm Spec}\,  \theta(g) =  \Delta_1(p)$ or
$ \Delta_2(p)$, or $\dim\theta = (p - 1)/2$ with  ${\rm Spec}\,\theta(g) =  \Delta_1(p) \setminus {1}$ or $ \Delta_2(p) \setminus {1}$. If $\ell = 2$,
then $H = PSL(2,p)$ and only $(p - 1)/2$ has to be taken for  $\dim\theta$;

$(7)$ $H = PSL(2,p)$ and $\dim\theta = p - 1$.

${\rm Spec} \,\theta(g)$ consists of all the non-trivial p-roots of $1$ except in cases $(5)$  and $(6)$. In case $(2)$
the eigenvalue $1$ does not occur for g belonging to one of the two unipotent conjugacy classes
of H.\end{propo}

Note that in the  cases $(1)-(7)$  $g$  non-trivially normalizes no  abelian subgroup of $G$. In many other cases $g\in N_G(A)\setminus C_G(A)$ for a suitable abelian subgroup
$A$ of $G$. The problem of obtaining a similar result for the $p$-power order elements of $G$ is only solved for classical groups \cite{DM08}. The cases with $G$ of exceptional Lie type is still open.  

Lemma \ref{hh1} is a useful tool for dealing with imprimitive linear groups.
The following fact of Clifford theory is essential for understanding the structure of primitive linear groups:   

\bl{cc5} {\rm (Clifford)} Let F be an \ac field and  $G\subset GL_n( F)$ 
be a tensor-indecomposable  primitive \ir subgroup. Then every non-central normal subgroup of G is \irt 

Moreover,   $G$ has an \ir subgroup that either is a central product of non-abelian quasisimple groups isomorphic to each other or 
  $n=r^{2m}$ for a prime r and some integer $m>0$, G contains an \ir subgroup E of symplectic type with $Z(G)\subset E$
and $G/E$ is isomorphic to a subgroup of $Sp_{2n}(r)$.  \el

We say that a subgroup $E\subset GL_n(F)$ is of symplectic type of $E$ is absolutely irreducible, $E'$ consists of scalar matrices and $n=r^m$ for some primes $r$ and $m>0$.
 
If $G\subset GL_n(F)$ is  tensor-decomposable then $G$ is   a subgroup isomorphic to a central product $GL_k(F)\circ GL_{n/k}(F)$ for some 
divisor $k$ of $n$, $k\neq 1,n$. In fact, the latter group occurs as the image of the tensor product  $GL_k(F)\otimes GL_{n/k}(F)$. \itf  $G$ has projective \ir \reps into $GL_k(F)$ and $ GL_{n/k}(F)$
and if $g_1,g_2$ are the images of $g$ (unique up to scalar multiples) in $GL_k(F)$ and $ GL_{n/k}(F)$, respectively, then $g=g_1\otimes g_2$ (the Kronecker product), again up to scalar.
Scalar multiples do not affect the degree of the minimal \po of a matrix, so one can ignore the multiples. By \cite{z05}, if $\ell=0$ then the spectra of  $g_1,g_2$   contain a scalar multiples of $\Delta_1\setminus 1$ or $\Delta_2\setminus 1$. In addition,   
$|\Delta_i\cdot \Delta_j|=p$ and 
$|(\Delta_i\setminus 1)\cdot (\Delta_j\setminus 1)|\geq p-1$   \cite[Lemma 2.13]{z05}. 

One observes that if $o(g)=p>2$ is a prime then $o(g_1),o(g_2)\in\{1,p\}$. If $|g_1|,|g_2|=p$ and $\deg g_1>(p-1)/2$ and $\deg g_2>(p-1)/2$ then one easily observes that  $\deg(g_1\otimes g_2)=p$. This hints that the following special case of Problem \ref{11} is of particular interest:

\begin{prob}\label{s11}  For an \ir finite quasisimple subgroups $G\subset GL_n(F)$ and an odd prime $p$ dividing $|G|$, when does the minimal \po degree $\deg g$ of a non-central  element $g\in G$ of order $p$  not exceed   $(p-1)/2$?\end{prob}

This problem is solved for $\ell=0$ in \cite{z05}. Specifically, in this case $G\cong SL_2(p)$ and $n=\deg g= (p-1)/2$.
   
If $\ell>0$ then there are other quasi-simple groups $G$ with $\deg g= (p-1)/2$ for some  non-central  element $g\in G$ of order $p$. An example is $M_{11}$ for $p=11$ and $\ell=3$. 

Lemma \ref{cc5} to  large extent  reduces the minimal \po problem (and some other linear group problems) to the cases where a minimal non-central subgroup of $G$ is either 
a central product of copies of a quasi-simple group or has an \ir normal subgroup of symplectic type. These cases are  of different nature. The second case is simpler.
 Instead of considering individual groups of this type, one can uniformly consider the normalizer $N$ of $E$ (added by the scalar non-zero matrices) in    $ GL_n(F)$.
The group $N$ is often referred as the extraspecial normalizer.
 
 If $g\in E\setminus Z(E)$ then $o(g)=r$ and $\deg g=r$. Let $g\in N\setminus E$. If $|g|$ is a $p$-power, $p\neq r$ then $g$ is contained in a finite $p$-solvable \ir subgroup of $N$,
and the situation is treated by a version of the Hall-Higman theorem. A more precise information on $\deg g $ can be extracted from \cite{Is73}, see also \cite[\S 3]{DZ}.  
The case with $p=r$ is more complex. For $g$ of odd prime power see \cite{BZ}. 

If $G$ has an \ir normal subgroup isomorphic to the direct product of simple groups, Problem \ref{is12} reduces to a more special case where $G$ contains an \ir simple non-abelian group.  Due to the importance of this special case, we single out the following special case of  Problems \ref{is12}:

  \begin{prob}\label{t42}  Let G be an almost simple group and F an algebraically closed field.  Determine \ir F-\reps $\phi$ of G  such that $\deg\phi(g)=o(g)$ for all $g\in G$.
\end{prob} 

J. Thompson \cite[p. 556]{Th1} suggests a more general version of Problem \ref{t42}: given a subgroup $H$ of $G$, and an \irr of $G$,  when is $\phi|_H=\rho_H^{reg}\oplus \si$ for some \rep $\si$ of $H$?    

\section{Eigenvalue 1 problems: Special cases}

\def\hw{highest weight }

\subsection{Simple algebraic groups}

In this section we discuss results related to Problem \ref{p5a}: determine unisingular \ir \reps of  simple algebraic groups. 

Note that  the center of a simple algebraic group $\GG$   is finite and $\GG/Z(\GG)$ is simple as an abstract group. Among the simple algebraic groups $\GG_1$ with $\GG_1/Z(\GG_1)\cong \GG/Z(\GG)$  there is  the largest one, 
in the sense  that $\GG_1/Z(\GG_1)\cong \GG/Z(\GG)$ implies  $\GG$ to be a quotient of $\GG_1$. Then  $\GG_1$ is called {\it universal} (or simply connected).  
One writes $\GG=\GG(K)$ to specify a "realization" of $\GG$ over a specific field $K$, which we always assume algebraically closed. The properties of $\GG$   depend on the characteristic of $K$,
whereas the choice of an algebraically closed field of fixed characteristic is rather  immaterial. Note however that if $K$ is the algebraical closure of a finite field then $\GG(K)$ is a locally finite group.
We usually omit mentioning of $K$ when discuss the structure of $\GG(K)$ or its representations.   

The universal simple algebraic groups over a fixed \acf $K$ are in 1-1 correspondence with the  simple Lie algebras over $\CC$. So they are labeled in the same way as the latter: $A_n,B_n,C_n,D_n,$ $n=1,2...$, $E_6,E_7$, $E_8$, $F_4$, $G_2$. The subscript is the {\it rank of the group}. Their \ir $K$-\reps are paramerized 
as those of Lie algebras over $\CC$ by strings of non-negative integers $(a_1\ld a_n)$, where $n$ is the rank. However, the dimensions (and other properties)
of $F$- and $\mathbb{C}$-\reps with the same parameter may be very different. 

Let $\GG$ be a connected algebraic group. 
Lemma \ref{w0w} tells us that $\GG$ is unisingular \ii every \irr of $\GG$ has weight 0. Simple unisingular algebraic groups are determined in  \cite{gt03}.

We now turn to Problem \ref{p5a}. This   first  arose in \cite{Z90} as a tool for 
study an analogous problem for \reps of finite groups of Lie type over the complex numbers. A rather obvious but useful  observation is the following:

\bl{w0w} Let $\GG$ be a connected reductive algebraic group (in particular, simple), and let $\phi$ be an \irr of  $\GG$. Then $\phi$ is unisingular \ii $\phi $ has weight $0$.
\el

Recall that every  connected reductive algebraic group has a maximal connected abelian subgroup $T$, say, consisting of semisimple elements called a {\it maximal torus of} $\GG$. The maximal tori of  $\GG$ are conjugate and every semisimple element of $G$ is contained in a maximal torus. The weights of $\phi$ are \ir constituents of $\phi|_T$, and the trivial constituent is called weight 0. 
(To be precise, weights are elements of $\Hom(T, K^\times)\cong \ZZ^n$, a $\ZZ$-lattice of some rank $n$, which is also called   the rank of $\GG$. So the weight 0 is the zero element of 
this lattice.)   
So if $\phi$ has weight 0 then every semisimple element of $\GG$ has \ei 1. By Lemma \ref{ssu}, $\phi$ is unisingular.
 The inverse holds too.  (Indeed, suppose the contrary. Let $\lam_1\ld \lam_m$ be the weights of $\phi$ and let $T_i=\{t\in T:\lam_i(t)=1\}$ be the kernel of $\lam_i$ for $i=1\ld m$. Then $T_i$
is a proper closed subgroup of $T$ (as $\lam_i\neq 0$). Then $S=T_1\cup...\cup T_m\neq T$, and hence $\lam_i(t)\neq 1$ for every $t\in T$, $t\notin S$.)

Note that the set  $\Hom(T, K^\times)\cong \ZZ^n$ has a certain partial ordering and every \irr $\phi$ of $\GG$ has a unique highest weights with respect to this ordering;  this is 
called the {\it highest weight of} $\phi$. This determines $\phi$ and is used for a parameterization of the \ir \reps of $\GG$. 
 
\med
Thus, Problem \ref{p5a} is equivalent to the following one:

 \begin{prob}\label{wz0}  Let $\GG$ be a connected reductive algebraic group. Determine the \ir \reps of $\GG$ that have weight $0$.
  \end{prob}

This problem has a satisfactory solution if $\ell=0$.    In fact in this case
 there is the following simple criterion for $\phi$ to have weight $0$:

\bl{w01} Let $\GG$ be a   simple algebraic group, and $\phi\in\Irr  \GG$ with \hw $\om$. 
Suppose that $\ell=0$. Then $\phi$ has weight $0$ \ii $\om$ lies in the root lattice. 
If $\GG$ is simply connected (or universal) then $\phi$ has weight $0$ \ii $\phi (Z(\GG))=\Id$.\el

Suppose that $\ell>0$. In this case one cannot expect to obtain a full description of the \ir \reps of every simple algebraic group  that have weight 0. The situation is untreatable for tensor-decomposable \reps with many tensor factors and groups $\GG$ of arbitrary rank. The solution is available for groups of  rank 1, which are $SL_2(K)$.  The case where $\phi$ is a tensor product of exactly two tensor-indecomposable terms is considered in \cite{BZ}. The case with three or more terms is hardly treatable (except certain special cases).
 It would be useful to obtain some sufficient conditions of  the weight 0 existence.

In \cite{Z90}  is discussed the uniform version of Problem \ref{wz0}:

 \begin{prob}\label{wz1} Determine the simple algebraic groups $\GG$ whose all \reps  have weight $0$.  \end{prob}

 Observe that $Z(\GG)=1$ for such groups. (If $Z(\GG)\neq 1$ then there exists \ir \reps non-trivial on $Z(\GG)$, 
and they are definitely not unisingular.) If $|Z(\GG_1)|$  equals the index of the root lattice in the weight lattice for the universal cover $\GG_1$ of $\GG$ then every 
  \irr of $\GG$ trivial on $Z(\GG)$ has weight 0. For instance, this is the case for $\GG=SO_{2n+1}(F)$ with $\ell$ odd. 
Problem \ref{wz1} is easy for tensor-indecomposable \ir representations. In \cite{gt03} Problem \ref{wz1} is solved in a more general context of simple groups of Lie type. 

 Finally we observe that Problem \ref{is16} is expected to have  a well describable solution. For instance, for $G=SL_n(F)$ we have: 

\begin{theo}\label{th3}{\rm \cite[Theorem 1.1]{zv23}}  Let $\GG=SL_n(F)$, where   $F$ be is algebraically closed field of characteristic $\ell$. Let $\phi $ be an \irr of $\GG$,  and let $g\in \GG$ be a rational element.   Then $\phi(g)$
has \ei $1$,  unless one of the following holds:

$(1)$   $\dim\phi=n$;

$(2)$ $(n,|g|,\dim\phi)\in\{(6,9,20),(6,15,20)  ,(8,15, 56),(10,45,120),\\(10,45,630), (14,45,364)\}$.   \end{theo}

Observe that elements $g$ in Theorem \ref{th3} are conjugate to those of $SL_n(\ell)$ so the result is valid for $SL_n(q)$ with $\ell|q$. 
There are similar results for $Sp_{2n}(F)$ and $SO_{2n+1}(F)$, see  \cite{zv23}.

\subsection{Groups of Lie type: \reps over a field of defining characteristic}
\def\hw{highest weight }

Many problems on finite linear groups over finite fields requires analysis of \ir \reps of 
algebraic groups. This is due to a fundamental theorem by Steinberg. Let $G$ be a quasi-simple finite group of Lie type in defining characteristic $\ell>0$. This means that $G$ is a subgroup of a simple algebraic group ${\mathbf G}$   such that $G=\{g\in {\mathbf G}: f(g)=g\}$ for a suitable surjective homomorphism $f:{\mathbf G}\ra {\mathbf G}$. This $f$ is called a {\it Frobenius} or {\it Steinberg} morphism. Let $F$ be an \acf of characteristic $\ell$.  

\begin{theo}\label{tst} (R. Steinberg) Every \ir $F$-\rep $\phi$ of $G$ extends to ${\mathbf G}$. \end{theo}


This theorem is highly useful as it reduces (to large extent) the \rep theory of $G$
to that of ${\mathbf G}$; the latter contains a lot of efficient tools to be used. For instance, the \ir $F$-\reps of $G$ are parameterized in terms of those of  ${\mathbf G}$.


\begin{corol}\label{tt8} Let $G$ be a finite group of Lie type in defining characteristic $\ell$ and let $\mathbf{G}$ be an algebraic group such that $G=C_{\mathbf{G}}( f)$, where $f$ is a Steinberg morphism of $\mathbf{G}$. Let $\phi$ be an absolutely \ir $\ell$-\rep of $G$. Then $\phi$ extends to a \rep $\tau$, say, of $\mathbf{G}$, and if $\tau$ is unisingular then so is $\phi$. \end{corol}

 The first statement of the lemma is Theorem \ref{tst}.  The second one follows from the fact that every semisimple conjugacy class of $G$ meets a maximal torus 
of   $\mathbf{G}$ (and Lemma \ref{ssu}).

As an illustration consider the following problems:

What can be said on an \ir $F$-\rep $\phi$ of $G$ if $\phi(G)$ contains a matrix with all \eis 
of \mult 1? of \mult at most 2? with all but one \ei \mult equal 1? with the \ei 1 \mult equal to 1? with no \ei 1? 

All these questions have analogs for algebraic groups ${\mathbf G}$. Moreover, these have 
a more uniform nature (in a sense) and methods to approach.

As mentioned above, ${\mathbf G}$ has an  abelian subgroup $T$ that contains a conjugate of every $\ell'$-element of ${\mathbf G}$.
(This is not true for the finite group $G$.)  This result allows to study $\ell'$-elements in $G$ uniformly by the study of $T$. 
 
For instance, if $\phi(G)$ has a $\ell'$-element with all \ei multiplicities 1 then, extending  $\phi$ to ${\mathbf G}$, one concludes that all \ir constituents of $\phi(T)$ occur with \mult 1.
The \ir \reps of $\GG$ with this properties are determined in \cite{Se} and \cite{SZ1}. Note that the \irr with at most one weight has \mult greater than 1 are determined in \cite{TZ2}.

Next we turn to Problem \ref{t41} (assuming that $\ell$ coincides with the defining characteristic of $G$). 
At the moment a full solution is available only for groups $GL_n(2)$ and $Sp_{2n}(2)$ and $\ell=2$.  
 
In this case the \ir \reps are in bijection with the strings $(a_1\ld a_{n-1})$
with entries $0,1 $. For instance, $(0\ld 0)$ is the label of the trivial \rep and $(1\ld 1)$
is the label of the maximum dimension \rep of degree $2^{n(n-1)/2}$.

\begin{theo}\label{th1} {\rm \cite{z18}}  Let $\phi=\phi(a_1\ld a_{n-1})$, $(0\leq a_1\ld a_{n-1}\leq 1)$, be a non-trivial \irr of $GL_n(2)$ over $\overline{\FF}_2$. Then $\phi$ is unisingular   \ii \\

\begin{center}
$\sum_{i=1}^{n-1} a_i\cdot i\geq  n$ and $\,\,\,\,\sum_{i=1}^{n-1} a_i\cdot (n-i)\geq   n.$\end{center}
\end{theo}

Observe that $\phi$ extends to a \rep $\tau$, say, of $\mathbf{G}=SL_n(F)$  with \hw $(a_1\ld a_{n-1})$, and $\tau$ is tensor-indecomposable. \itf $\tau$ has weight 0 \ii 
$\sum_{i=1}^{n-1} a_i\cdot i\equiv 0 \pmod{n}$. So there are many unisingular \ir \reps of $G$ that do not extend to   unisingular \ir \reps of  $\mathbf{G}$. 

The result for $G=Sp_{2n}(2)$ is more complex to state; however, for prime power order elements we have:

\begin{theo}\label{pp1} {\rm \cite[Theorem 1.4]{z22}}  Let $G=Sp_n(2)$ and let   $g\in G$ be an element of odd prime power order. Let $\phi=\phi(a_1\ld a_{n-1})$, $(0\leq a_1\ld a_{n}\leq 1)$, be a non-trivial \irr of  G over $\overline{\FF}_2$. 
If $ a_1+\cdots + a_{n}>1$ then $1$ is an \ei of $\phi(g)$. 
\end{theo}

The method used in {\rm \cite{z18}} allows to obtain  some sufficient conditions for an \irr of finite simple group of Lie type to be unisingular. The following result is 
proved in \cite{cz}.

\begin{theo} \label{s21}
Let $G$ be a finite simple group of Lie type $A_n(q)$, $C_n(q)$ or $D_n^\pm (q)$, q a prime power,
$\mathbf{G}$ the respective universal simple algebraic group and let $\om_1\ld \om_n$ be the fundamental weights  of $\mathbf{G}$. Let $V$ be an \ir $\mathbf{G}$-module. Suppose that the \f conditions hold:

$(1)$ $\mathbf{G}$ is of type $A_n$ and there  are natural numbers $m_1,m_2$ and $i\in\{1\ld n\}$ such that $m_1(q-1)\om_i+\om_1+\om_n$ and $m_2(q-1)\om_i$ are weights of $V$.

$(2)$ $\mathbf{G}$ is of type $C_n$ $(n>1)$ or $D_n$ $(n>3)$ and there are natural numbers $m_1,m_2,m_3$ such that  
$m_1(q+1)\om_1, m_2(q-1)\om_1$ and $m_3(q-1)\om_1+\om_2$ are weights of V.

\noindent Then  every semisimple element $g\in G$ has \ei $1$ on V.\end{theo}


The approach used in the proof of this result requires  a further improvement to deal with the unitary groups. 
 
\subsection{Cross-characteristic representations}

Not much is known on Problem \ref{t41} for cross-characteristic representations of simple groups $G$ of Lie type. Recall that this term means that the  characteristic $\ell $ of the representation field is distinct from the defining characteristic of a group of Lie type. If the $\ell$-decomposition matrix of $G$  is known one can in principle obtain an explicit solution   to Problem \ref{t41},
provided it has been solved for the \reps over the complex numbers. However, even for the latter  case the problem remains difficult, and the results are far to be complete.

Some information can be extracted from \cite{TZ08,TZ22} and \cite{TZ25}. In particular, it is shown that the $\ell$-Brauer characters of certain simple groups of Lie type 
are constant on the regular elements of some maximal tori. See also \cite[Remark 2.3.10]{MT}.

There are some results on the fixed point space dimension     of some elements  in an \ir linear group $G\subset GL(V)$ \cite{sh,LSh}. 
In \cite{LSh} it is proved that if $|G/R|>120$, where $R$ is a maximal solvable normal subgroup of $G$, then  $\dim C_V(g)\geq \dim V/6$ for some $g\in G$. If $G$ is simple then 
 \cite[Theorem 1.3]{GM} states that $\dim C_V(g)\leq \dim V/3$ for some $g\in G$.


\section{Linear groups and group representations over the complex numbers}

The problem of determining 
unisingular \ir \reps of   finite simple groups over the complex numbers is in general open. According to our knowledge,
 a complete solution to Problem \ref{t41} is available only for the class of symmetric and alternating groups \cite{apv, ap2}. In particular,  if  $G=S_n$, $\dim\phi>1$ and  $n>10$ then $\deg\phi(g)\geq |g|-1$ with equality only if $\dim\phi=n-1$ and $g\in G$ is a full cycle. See also \cite{st1} where the author obtains a solution to    Problem \ref{is12} for representations of symmetric groups over the complex numbers. See also \cite{pv3}. 

For simple groups $G$ of Lie type the problem was first discussed in \cite{Z90}. Let $\ell$ be the defining characteristic of $G$. The results obtained there were based on the fact that
there are groups $G$ whose all $\ell$-modular \reps are unisingular, equivalently, have weight 0. (A full list of such groups is obtained in \cite{gt03}.) Using reduction of \reps of $G$ over the complex numbers to characteristic $\ell$, one concludes that all $\ell'$-elements  of $G$ have \ei 1.    

So for such groups the problem reduces to analyzing the  $\ell$-singular elements. 

The minimal polynomials and \eis of $\ell$-elements of quasisimple groups of Lie type are studied in  Proposition \ref{z90} for elements of order $\ell$ (for cross-characteristic representations),
elements of $\ell$-power order are treated in \cite{DM08}. See also  \cite{z05} and   \cite{z08}.
 
The $\ell$-singular elements behaves particularly well if $\phi$ is the Steinberg \rept (The latter is the only 
\irr of a simple group of Lie type of $\ell$-defect 0, equivalently, of degree $|G|_\ell$, the order of a $\ell$-Sylow subgroup of $G$; so this property can be taken for a definition of the Steinberg \rep when $G$ is simple.) 
Then it was deduced that for   $\ell$ odd the Steinberg \rep of $G$ is unisingular. The case of $\ell$ even remained open until recent years: some cases has been resolved in     \cite{cz}, the remaining cases were considered in \cite[Theorem 1.1]{z24}.
The final result accumulated the above contributions is: 

\begin{theo} \label{st22} Let $G$ be a simple group of Lie type and $\si$  the Steinberg \rep of G. Then $\si$ is unisingular except for the case where $G=PSL_2(q)$ with q even. \end{theo}
 
We mention an interpretation of the unisingularity problem (for $F$ of characteristic 0)  in terms of \rep theory. This is based on the Frobenius reciprocity theorem applied to the cyclic group $C=\lan g \ran$. 

\bl{fr1} An \irr $\rho$ of G is a constituent of the induced \rep $1_C^G$ \ii
$\rho(g)$ has \ei $1$. Consequently, $\phi$ is unisingular \ii $\phi$ is a constituent
of $1_C^G$ for every cyclic subgroup C of G. \el

There is a natural problem which is in a  sense dual to the unisingularity problem:

\bpr{is27}   Determine the set  of  elements $g\in G$ that has \ei $1$ in every \irr of G.\epr

This is equivalent to 

\bpr{is28}  Determine cyclic subgroups $C$ of $ G$ such that every \irr of G is a constituent of the induced \rep $1_C^G$.\epr

In this form Problem \ref{is28} is discussed in \cite{hstz} and \cite{pv3}.
This is related with the problem of decomposing the permutation characters associated with the conjugation action of a group $G$ on its conjugacy classes;
this  is equivalent of decomposing the induced characters $1_{C_G(C)}^G$. 
This is discussed in \cite{hz06} and  \cite{hstz} and elsewhere. 


\section{Minimal groups and abelian-by-cyclic groups}

A general principle to approach a new problem is to start from a minimal non-trivial special case. In this context we consider non-abelian solvable groups with only one non-trivial proper subgroup.
It is clear that such a group $G$ has a normal elementary abelian  $r$-subgroup $A$, say, for some prime $r$ and $G/A$ is cyclic of prime order $p\neq r$. (So $G$ is a semidirect product of $A$ and $G/A$.)
Moreover, $A$  has no proper non-trivial subgroup that is normal in $G$. Recall that an elementary abelian $r$-group can be viewed as an additive group of a vector space $V$ over 
$\FF_r$ whose dimension is equal to the rank of $A$ and the conjugation action of $G$ on $A$ turns $V$ to an $\FF_r G$-module, in fact, an $\FF_r(G/A)$-module. The condition that 
$A$  has no proper non-trivial subgroup   of $G$ is interpreted by saying that $A$ is an irreducible    $\FF_r(G/A)$-module. This interpretations is rather useful. Indeed, let $H=G/A=\lan h\ran$, that is, $h$ is a generator of $H$. Then  $h$ acts irreducibly on $V$ \ii $n=\dim V$ is the minimal number $d$ such that $p|(r^d-1)$. This number is called the order (or the multiplicative order) of $r$ modulo $p$ and denoted usually by ${\rm ord}_p(r)$. Note that   ${\rm ord}_p(r)$ can be also defined as the minimal number $d$ such that $\FF^\times _{p^d}$ contains an elements of order $p$,
or as $\FF_r(\zeta):\FF_r$, where $\zeta$ is a primitive $p$-roots of unity in an algebraically closed field of characteristic $r$. By Fermat's little theorem, $p|(r^{p-1}-1)$, so ${\rm ord}_p(r) \leq p-1$. Therefore,  $\FF_r(\zeta)$ is a subfield of $\FF_{r^{p-1}}$, and hence ${\rm ord}_p(r)$  is a divisor of $p-1$.

One observes that the group $G$ just constructed depends on two parameters $r,p$ as ${\rm ord}_p(r)$ is a well defined number. We denote this group by $G_{r,p}$. (This is a special case of a Frobenius group \cite[\S 25]{Fei}.)
Note that the group $G_{r,p}$ can be defined in terms of a finite field as follows. Let $d={\rm ord}_p(r)$, $q=r^d$, $H$ a subgroup of $ \FF^\times _q$ of order $p$  and $A=\FF^+ _q$ be the additive group of $\FF  _q$.
Then $G:=HA$, the semidirect product, is isomorphic $G_{r,p}$. Note that the rank of $A$ equals  ${\rm ord}_p(r)$.  

\med
Example.   ${\rm ord}_3( 5)=4$,  ${\rm ord}_3( 11)=5$   and ${\rm ord}_3(23)=12$.

\bpr{is17}   For $r>2$ determine unisingular \reps of $G_{r,p}$ (over $\CC$ or any field of characteristic $\ell\neq r$).\epr

Note that Problem \ref{is17} is a special case of Problems \ref{12} and \ref{is13}. We assume $r>2$ as every faithful \irr of $G_{2,p}$ is unisingular.

One can expect that Problem \ref{is17} is easy as the group itself is of very simple nature. All elements of $G$ are of order $p,r$ or 1, and all \ir \reps over an algebraically closed field are 
of degrees $p$ or 1, those of degree $p$ are faithful and lift to characteristic 0. So it suffices to deal with \reps over the complex numbers.    In addition,  all elements of order $p$ have \ei 1. 
So Problem \ref{is17} in fact asks 
whether (or when)   every element of order $r$ has \ei 1. They all lie in an abelian normal subgroup $A$.  If $\ell\neq r$ then $\phi(A)$ is diagonalizable in any $\ell$-\rep of $G$.
We state the \f easy lemma: 

\bl{aa1} Let $A\subset GL_n(F)$ be an elementary abelian $r$-group of diagonal matrices.  Let $K_i$ be the subgroup of $A$ of elements that have $1$ at the $(i,i)$ diagonal position.
Then $A$ is unisingular \ii $A=\cup K_i$.  
\el

 Next suppose that $G$ is a group of monomial matrices in $GL_n(F)$ (this is exactly the normalizer in $GL_n(F)$ of the group $D$ of diagonal matrices). Set  $A=G\cap D$. 
Then $G$ permutes the homogeneous components of $A$, and in particular, $G$ acts by permutations on the set of \ir constituents of $A$, and also permutes their kernels.
 If $G$ is \ir  then they form a single $G$-orbit. If $K$  is one of them   then $\cap_{ g\in G}\, gKg\up =1$. In addition, $A$, and hence $G$,  is unisingular \ii  $A=\cup_{g\in G}\, gKg\up$.
So the question when $A=\cup_{h\in H} hK_1h\up$
is a special case of Problem \ref{is19}.  

Viewing $A$ as a vector space $V$, say, over $\FF_r$ and $K$ as a subspace $W$ of codimension 1, we observe that  $A=\cup_{h\in H} gKg\up$ can be rewritten as $V=\cup_{h\in H} hW.$ We often record this as $V=HW$.
 
Next we turn to groups $G=G_{r,p}=AH$, where $A$ is an elementary abelian $r$-group, $H$ is of prime order $p$ and $A$ is the only non-trivial proper normal subgroup of $G$. 
This group is abelian-by-cyclic. \itf every \irr  of $G$ (over an algebraically closed field $F$)  is either one-dimensional or of degree $p$, 
in fact, is  an induced \rep  $\lam^G$, where $\lam$ is a non-trivial  \ir (hence one-dimensional) \rep of $A$.  By Clifford's theorem, the restriction of  $\lam^G$ to $A$ is a direct sum of 
$p$ distinct \reps $\lam^h$, where $\lam^h(a)=\lam(hah\up)$, $h\in H$. Furthermore, if $\lam,\mu\in \Irr A $ then $\lam^G$ and $\mu^G$ are equivalent \ii the restrictions $\lam^G|_A$
and $\mu^G|_A$ have no common \ir constituents. Therefore, the \ir \reps of $G$ of degree $p$ are in bijection with the non-trivial $H$-orbits on $\Irr A$.  

Representations $\tau,\si$ of a group $G$ are said to be quasi-equivalent if there exists an \au $\al$ of $G$ such that $\tau$ is equivalent to $\si^\al$, where $ \si^\al(g)=\si(\al(g))$. 
We show that the faithful \ir \reps of $G_{r,p}$   are quasi-equivalent (so they behave uniformly). In particular, either all such \reps are  unisingular or all they are not unisingular. 

\bl{rq2} Faithful \ir \reps of $G=G_{r,p}$ are quasi-equivalent. \el

\bp We have $G=A.H$, where $A$ is an elementary abelian normal $r$-subgroup of $G$,   $|H|=p$ and $H$ acts on $A$ as a subgroup of $GL(A)$. 
Let $d$ be the rank of $A$. We have seen above that $d={\rm ord}_r(p)$, so $H\subset GL(A)\cong GL_d(r)$.  Let 
$C=C_{GL_d(r)}(H)$ so $H\subseteq C$. Then $C$ normalizes $G$. So the conjugation by elements of $C$ yield automorphisms of $G$.  

As $H$ is \ir in $GL_d(r)$, by Schur's lemma,  $\FF_rH$ (the $\FF_r$-enveloping algebra of  $H$ in $\End V$) is a field isomorphic to $F_{r^d}$, 
and hence $C\cong F^\times_{r^d}$.  Moreover, $C$ acts transitively on $V\setminus \{0\}$. Note that    $ \Irr A\cong  \Hom(A,F^\times)$ is the dual group of $A$.  
Then $C$ acts transitively on the non-trivial elements of $\Irr A$, so if $\lam,\mu\in\Irr A$ are non-trivial then $\mu=\lam^c$ for some $c\in C$. Therefore, if $\mu^G=(\lam ^c)^G=(\lam^G)^c$
as required. \enp
  
The \f is now obvious: 

\bl{qe1} Suppose that \reps $\tau,\si$ of a group G are quasi-equivalent. Then  $\tau$ is unisingular \ii so is $\si$.
\el

This implies that if $H\subset GL(V)$ is an \ir cyclic subgroup and $V=HW$ for a subspace $W$ of codimension 1 in $V$ then $V=HW'$ for every such subspace $W'$.  
 This fact   has an essencial computational significance. Indeed, in order to check by computer the equality  $V=HW$ for every subspace $W$ of codimension 1, it suffices to 
choose $W$ in random.  In particular, Problem \ref{is17} can be stated in the \f form:     
 
\bpr{is17a} 
 For $r>2$ determine the pairs $(r,p)$ such that $G_{r,p}$ has a unisingular \irr (over $\CC$ or any field of characteristic $\ell\neq r$).\epr

 Recall that $d={\rm ord}_r(p)$ divides $p-1$.  As we have already seen,  this is equivalent to a vector space covering problem:

\bpr{is17b} Determine the pairs  $r,p$ of primes,  $r>2$, such that if  $H\subset GL_d(r)=GL(V)$ be an \ir subgroup of order p then $V= HW$ for some subspace $W$ of codimension $1$ in V.  \epr
 
Thus, one can approach Problem \ref{is17} using either matrix theory or vector space theory. In fact, the matrix approach can be replaced by the \rep theory approach. Both of them are not 
very efficient but each of them has some advantage in some special cases.  Related problem is under discussion in \cite{cer} and \cite{agt}. 

  Problems \ref{is17a}  and \ref{is17b} seems to be difficult already for $r=3$, see \cite{z24}. The special case where $d=p-1$ is settled in  {\rm \cite[Lemma 1.5]{z24}} as follows:

\begin{theo}\label{th4}  Suppose that ${\rm ord}_r(p)=p-1$. Then $G_{r,p}$ has no non-trivial unisingular representation. \end{theo} 

\med
For very small $p$ there are character tables of $G_{r,p}$, so one can easily answer the question whenever the character table is available.   However, the situation looks rather obscure in general.  
  
Keep $d={\rm ord}_r(p)$. The case with $d=1$ is trivial, as  $A$ is cyclic. One easily observes that in this case the only unisingular \rep is the trivial one.  
The case $d=2$ is easy too. In this case $p|(r+1)$; here $p\neq r+1$ as $p,r$ are odd, and hence $p<r$. In this case $G_{r,p}$ has no unisingular 
non-trivial  \irr due to the following easy observation:

\bl{eo1} Let $p<r$. Then $A\neq \cup_{h\in H}hKh\up$.\el

\bp We have $|A|=r^d$ and $|K|=r^{d-1}$. Then  $|\cup_{h\in H}hKh\up|\leq p|K|<r|K|=|A|$.\enp

Note that the inequality  $p<r$ can be improved to $p<3(r+1)/2$, see  \cite[Lemmas 2.1 and  2.2]{Be}.  
For our case, where $p$ is a prime, we have $p\neq 3(r+1)/2$.

\bpr{is17c} Investigate Problem $\ref{is17b}$ for the special case   with   $d=(p-1)/2$.\epr  

It has been verified by computer that the \ir \reps of degree greater that 1 of  the groups  $G_{3,11}$ and $ G_{3,23}$ are  unisigular, see \cite[Lemma 6.3]{z24}. 
   Computer power seems to be  insufficient to check the case with  $G_{3,47}$.

In view of Proposition \ref{es3}, the group $G=G_{r,p}$ can be realized as a transitive subgroup of the symmetric group $S_{rp}$ and  $G_{r,p}$  has no non-trivial
unisingular \ir \rep \ii $G$ has a fixed point free element of order $r$. In \cite{cer} the authors state the following conjecture which they refer to as Isbell's conjecture:

\med
{\bf Conjecture} Let $H$ be a transitive subgroup of $S_n$, where $n=r^mp$, $r,p$ are primes. Then there exists a function $f=f(r,p)$ such that $H$ has a fixed point free $r$-element whenever $m>f$. 

\med
If $H=G_{r,p}$ then the Sylow $r$-subgroup of $H$ is elementary abelian, and $d|(p-1)$. So Problem \ref{is17c} does not interact with Isbell's conjecture.

In view of Theorem \ref{th4}  one can ask, given a prime $r$, whether there are infinitely many primes
$p$ such that ${\rm ord}_r(p)=p-1$. This is known as Artin's primitive root problem which remains open since 1929. 

\bpr{ar17} Let $r>2$. Is it true that there are infinitely many primes p such that  some faithful \ir \rep  of $G_{r,p}$ is not unisingular?
 \epr

We expect \cite[Conjecture 3]{z24} that the answer is positive.
 
\bigskip

Conclusion. The research area related with the \ei 1 problems of elements in group \reps
is at the initial stage of study, and   promising for further progress.

 \bigskip

{\it Conflict interest statement} There is no conflict of interest with any person or organization.


 \end{document}